\documentclass[letterpaper, 10 pt, conference]{ieeeconf}  

\IEEEoverridecommandlockouts                              

\usepackage{amsmath,mathtools}
\usepackage{amssymb}
\usepackage{times}
\usepackage{multirow}
\usepackage{graphicx}
\usepackage[dvipsnames,table,xcdraw]{xcolor}
\usepackage{tikz}
\usepackage{subcaption}
\usepackage[normalem]{ulem}
\usepackage{algorithmic,algorithm}
\usepackage[noadjust]{cite}
\usepackage[hyphens]{url}
\usepackage{hyperref}
\usepackage[hyphenbreaks]{breakurl}
\usetikzlibrary{arrows}
\usetikzlibrary{calc}

\newcommand{\Xs}{\mathcal{X}}

\newcommand{\Sem}{\textbf{\emph{s}}}
\newcommand{\Aem}{\textbf{\emph{a}}}
\newcommand{\Iem}{\textbf{\emph{i}}}
\newcommand{\Rem}{\textbf{\emph{r}}}
\newcommand{\Dem}{\textbf{\emph{d}}}

\newcommand*\diff{\mathop{}\!\mathrm{d}}

\newcommand{\yonote}[1]{\textcolor{red}{[YO: #1]}}

\def\Re{\mathbb{R}}

\def\argmin{\mathop{\text{\rm arg\,min}}}

\def\Sec#1{Sec.~\ref{#1}}

\def\notes#1{\marginpar{\tiny #1}\typeout{Notes!
Notes!
Notes!
}}
\renewcommand{\notes}[1]{\typeout{notes!}}

\def\Re{\field{R}}

\def\Sec#1{Sec.~\ref{#1}}

\def\Sec#1{Sec~\ref{#1}}

\def\E{{\sf E}}

\def\Sec#1{Sec.~\ref{#1}}

\newtheorem{assumption}{Assumption}
\newtheorem{definition}{Definition}

\newtheorem{remark}{Remark}
\newtheorem{proposition}{Proposition}

\def\beq{\begin{eqnarray}} 
\def\bc{\begin{center}} 
\def\be{\begin{enumerate}}
\def\bi{\begin{itemize}} 
\def\bs{\begin{small}}
\def\bS{\begin{slide}}
\def\ec{\end{center}} 
\def\ee{\end{enumerate}}
\def\ei{\end{itemize}}
\def\es{\end{small}}
\def\eS{\end{slide}}
\def\eeq{\end{eqnarray}}

\newcommand{\newP}[1]{\medskip\noindent{\bf #1:}}

\newcommand{\ud}{\,\mathrm{d}}

\def\Re{\mathbb{R}}
\def\E{{\sf E}}

\def\argmin{\mathop{\text{\rm arg\,min}}}

\def\Sec#1{Sec.~\ref{#1}}

\def\Prop#1{Prop.~\ref{#1}}

\renewcommand{\Re}{\mathbb{R}}

\title{\LARGE \bf
Modeling Presymptomatic Spread in Epidemics via Mean-Field Games
}

\author{S. Yagiz Olmez, Shubham Aggarwal, Jin Won Kim, Erik Miehling
  \\[5pt] Tamer Ba{\c s}ar, Matthew West, and Prashant G. Mehta\\[5pt]
\thanks{Research supported in part by the C3.ai Digital Transformation Institute sponsored by C3.ai Inc. and the Microsoft Corporation and in part by the National Science Foundation grants NSF-ECCS 20-32321 and NSF-CMMI 1761622.}
\thanks{	S. Y. Olmez, S.Aggarwal, J. W. Kim, and P. G. Mehta are with the Coordinated Science Laboratory and the Department of Mechanical Science and Engineering at the University of Illinois at Urbana-Champaign (UIUC); E. Miehling and T. Ba{\c s}ar are with the Coordinated Science Laboratory and the Department of Electrical and Computer Engineering at UIUC; M. West is with the Department of Mechanical Science and Engineering at UIUC; Corresponding email: mehtapg@illinois.edu.}
Coordinated Science Lab, University of Illinois at Urbana-Champaign, Urbana, IL\\%
{\tt\small \{solmez2,sa57,jkim684,miehling,basar1,mwest,mehtapg\}@illinois.edu}
}

\begin{document}

\tikzstyle{circ} = [draw,circle,fill = white!20,minimum width = 3pt, inner sep  = 5pt]
\tikzstyle{line} = [draw, -latex']

\maketitle
\thispagestyle{empty}
\pagestyle{empty}
	\begin{abstract}
This paper is concerned with developing mean-field
game models for the evolution of epidemics. Specifically, an agent's decision -- to be socially active in the midst of an epidemic -- is modeled as a mean-field game with health-related costs and activity-related rewards.  By
considering the fully and partially observed versions of this problem, the role of information in guiding an agent's rational decision is highlighted. The main contributions of the paper are to derive the equations for the mean-field game in both fully and partially observed settings of the problem, to present a complete analysis of the fully observed case, and to present some analytical results for the partially observed case. 

	\end{abstract}

\section{INTRODUCTION}

It is often asked what is so hard about modeling the spread of COVID-19?
Although modeling proved to be invaluable during the early stages of
the pandemic (March 2020), primarily by persuading the reluctant
politicians to adopt the harsh lockdown measures~\cite{nigel2020illinois}, subsequent 
evolution of the pandemic has shone a harsh spotlight on the
simplistic SIR models~\cite{nyt2021india}.  This
in turn has spurred much work on its enhancements by including
additional compartments~\cite{giordano2020sidarthe,fernandez2020estimating,peirlinck2020visualizing,anastassopoulou2020data,chen2020time,ku2020epidemiological,bastos2020modeling}, considering the effects of
heterogeneity~\cite{tkachenko2021herd,britton2021hetero}, and modifying the mean-field interaction terms~\cite{taghvaei2020fractional}.  

There are three aspects to modeling the pandemic: (i) evolution of virus
in a single agent who has been infected; (ii) social behavior of a 
single agent up to the time that she is infected; and (iii) the net
effect (mean-field) due to the social behavior of the population.  Our understanding of
the aspect (i) has improved by leaps and bounds based on COVID-19
data~\cite{mettler2020data,forecasting2020,numbers2020,byrne2020,olmez2021data}.
Aspects (ii) and (iii) have proved to be
nearly intractable and this has prompted resorting to large-scale
agent-based simulations with its many assumptions and
paramaters~\cite{imperial2020agent,gopalan2020reliable,chang2020modelling,liu2020machine,tkachenko2021,ozik2021argonne}.

The mathematical modeling and analysis of the aspects (ii) and (iii)
is an extremely complex problem of immense societal importance.  This
paper is a modest effort in this direction based on the mean-field
game (MFG)
formalism.  Specifically, the agent's social behavior (aspect (ii)) is
modeled as an optimal control problem based on health-related costs
and activity-related rewards.  The collective effect of the population (aspect (iii)) is modeled in terms of two mean-field
processes $\beta$ and $\alpha$.  The instantaneous $\beta_t$ in
particular models the effect of the active infected agents.  



The main contributions of the paper are to derive the
equations for MFG in both fully and partially observed
settings of the problem, to present complete analysis of the fully
observed case, and to present some analytical results for the
partially observed case. Specifically, the MFG model is used to obtain
the following two conclusions:

\medskip

\noindent
\textbf{Rationality of a single agent:} A susceptible agent will
choose to be active if and only if the reward outweighs the risk
(formalized by deriving a certain critical value
$\beta^{\text{crit}}$ for $\beta$).  In
contrast, a rational infected agent who also perfectly knows her
epidemiological status will choose to self-isolate (quarantine).  The
latter guided the testing regimens deployed by the
U.S. universities during the fall and spring of 2020~\cite{nyt2021davis,nyt2021beyond}.
Unfortunately, the observed behavior of partying U.S. undergraduates
served to highlight the limitations of rationality~\cite{fall2020covid19,nyt2021notgreat}.      

\medskip

\noindent
\textbf{Imperfect information:} Under imperfect information -- when an
agent does not have perfect belief regarding her own epidemiological
status -- an infected agent may behave as a susceptible agent.  Such a
behavior may in turn drive the epidemic which is consistent with the
reported role of
presymptomatic and asymptomatic population in the COVID-19 spread
\cite{rivett2020screening,buitrago2020role,bender2021analysis,wei2020presymptomatic}.


\medskip
The importance of individual decision-making affecting the evolution
of the epidemic was recognized early on with contributions on MFG
modeling of epidemics appearing
in~\cite{hubert2020incentives,lee2020controlling,elie2020contact,cho2020mean,aurell2020optimal,doncel2020mean,aurell2021finite,tembine2020covid}. Related
to our work are~\cite{elie2020contact} where an agent's
decision variable is her rate of contact with others, and~\cite{aurell2020optimal,aurell2021finite} where an agent strives to
follow a prescribed rate of contact based on government guidelines.
The novelty of our paper comes from partially observed settings
and differences in cost structures which are helpful to model and analyze the presymptomatic spread of
epidemic.


The remainder of this paper is organized as follows: The problem
formulation appears in \Sec{sec:problem}.  Its solution for the fully
observed and partially observed cases is described in \Sec{sec:fully}
and \Sec{sec:partially}, respectively.  The proofs appear in the
Appendix.

\def\b{\beta}
\def\v{\alpha}
\def\f{\bar{\phi}}

\begin{figure*}[!t]
	\centering
	\begin{tabular}{ccc}
		\begin{subfigure}[b]{0.28\textwidth}
		\centering
 		\begin{tikzpicture}[node distance = 2cm,auto]
        \node[circ] (S) {$\Sem$};
        \node[circ, above of = S] (A) {$\Aem$};
        \node[circ, right of = A] (I) {$\Iem$};
        \node[circ, right of = I] (D) {$\Dem$};
        \node[circ, below of = I] (R) {$\Rem$};
        \path[line] (S) -- node[midway,left] {$\lambda^{\text{\tiny SA}} \beta_t  U_t$} (A);
        \path[line] (A) -- node[midway,above] {$\lambda^{\text{\tiny AI}}$} (I);
        \path[line] (I) -- node[midway,above] {$\lambda^{\text{\tiny ID}}$} (D);
        \path[line] (I) -- node[midway,right] {$\lambda^{\text{\tiny IR}}$} (R);
        \node at (-1.7,2.0) {presymptomatic};
        \node at (-1.3,0) {susceptible};
        \node at (3.2,0) {recovered};
        \node at (4.9,2) {dead};
        \end{tikzpicture}
        \caption{}
 		\end{subfigure}
		&$\quad\quad\quad\quad\quad\quad$&
		\begin{subfigure}[b]{0.38\textwidth}
		\centering
		\begin{tabular}{|c|c|c|c|}
		\hline
		\multirow{3}{*}{\textbf{state}} & \textbf{health} & \textbf{alturistic} & \textbf{activity} \\
		& \textbf{cost} & \textbf{cost} & \textbf{reward} \\ \cline{2-4}
		& $c^{\text h}(x)$ & $c^{\text a}(x)$ & $r(x;\v)$ \\
		\hline
		$\Sem$ & $0$ & $0$ & $\v$ \\
		$\Aem$ & $0$ & $1$ & $\v$\\
		$\Iem$ & $1$ & $1$ & $\v$ \\
		\hline
		\end{tabular}
		\caption{}
		\end{subfigure}
	\end{tabular}
\caption{(a) Epidemiological states and the transition graph. (b) Cost
  function. }
\label{fig:model}
\end{figure*}
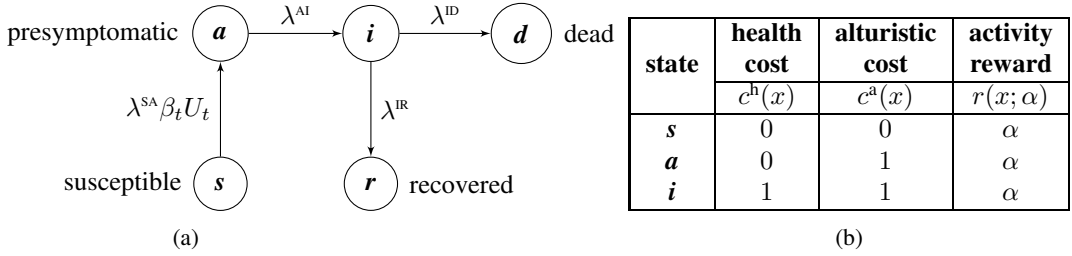

\section{Problem formulation: Modeling}\label{sec:problem}

\subsection{Model for a single agent}
\label{sec:model_single}
\newP{Dynamics} For a single agent, the epidemiological state is
modeled as a Markov process and denoted by $X:=\{X_t\in\Xs: t\geq 0\}$
where the state-space is $\Xs := \{\Sem,\Aem,\Iem,\Rem,\Dem\}$.
Figure~\ref{fig:model} depicts the transition graph and includes a
description of the
epidemiological meaning of each of the five states.
There are two types of infected states: (i) presymptomatic, denoted as
$\Aem$; and (ii) symptomatic, denoted as $\Iem$.  In either of these
states, an agent is 
infectious, i.e., able to infect other agents. The modeling distinction is
that, in partially observed
settings, a presymptomatic agent may not {\em know} her true state but a symptomatic agent does. Note that, for the present paper, the transition graph does not include an edge from $\Aem$ to $\Rem$ and thus one may interpret asymptomatic as presymptomatic whereby the agent discovers that she is infected upon appearance of the
symptoms. Broadly, there are two types of transitions:

\medskip

\noindent \textbf{1)} On the subset $\{\Aem,\Iem,\Rem,\Dem\}$, the
transition rate depends {\em only} upon the agent attribute $\theta$,
which here represents the age of the agent.  For example, an older
infected agent may risk a longer recovery time (smaller
$\lambda^{\text{\tiny IR}})$ than a younger agent.

\medskip

\noindent \textbf{2)} The transition from $\Sem \to \Aem$ depends upon three factors:
  (i) the intrinsic infectivity of the virus; (ii) the agent behavior
  (level of social activity); and (iii) the behavior of the
  infected agents in the population. The following equation is used
  to model the effect of these three factors:
 \[
 \text{rate}[\Sem \to \Aem] = \lambda^{\text{\tiny SA}}\ \b_t\,  U_t 
 \]
 where $U=\{U_t\in[0,1]:t\geq 0\}$ is the agent's activity {level}.  For
 $U_t=0$ (resp., $U_t=1$) the agent is completely isolated
 (resp., active) at time $t$.  The rate also depends upon the
 process $\b:=\{\b_t\in[0,1]: t\geq 0\}$ which is used to model the
 net effect of the behavior (activity) of infected agents {in the population}.  Its
 model is introduced in \Sec{sec:pomfg_model} together with the population model. Finally, the
 parameter $\lambda^{\text{\tiny SA}}$ models the rate of spread of
 the virus in a given population in the absence of any mitigation measures.


\medskip

\noindent \textbf{Optimal control objective:}  In the following,
$\v:=\{\v_t\in[0,1] : t\geq 0\}$ and $\b$
are given deterministic
processes. The control objective for a single agent is to chose her
activity $U$ to minimize 
\[
{\sf J}(U;\b,\v)  = {\sf E} \left( \int_0^T
  e^{-\gamma t} c(X_t,U_t;\v_t) \ud t + e^{-\gamma T} \phi(X_T) \right)
\]
where $T=T(\omega)=\inf\{t>0: X_t(\omega) \in \{\Dem,\Rem\}\}$ is the
random stopping time when the agent either recovers ($X_T=\Rem$) or the
agent dies ($X_T=\Dem$); by convention, $\inf \emptyset = \infty$.  The cost
function is of the following form:
%
\begin{align*}
c(x,u;\v_t) &= c^{\text{h}} (x)+c^{\text{a}} (x) u -r(x;\v_t)u
\end{align*}
where models for the health related costs $c^{\text{h}}$ and 
$c^{\text{a}}$, the activity related reward $r(x;\cdot)$, and the
terminal cost 
$\phi(x)$ are tabulated in the Table included as part of Fig.~\ref{fig:model} (b).

\medskip

\noindent \textbf{Information structure:} There are two settings of
the problem: (i) the fully observed case; and (ii) the partially
observed case. 
In the partially observed setting, the observation process 
$Y:=\{Y_t\in\{0,1\}^2:t\geq 0\}$ is defined according to 
\[
Y_t = \begin{bmatrix} 1_{[X_t=\Iem]} & 1_{[X_t=\Dem]} \end{bmatrix}
\]  



\subsection{Model for the mean-field}
\label{sec:pomfg_model}

To fully specify the problem, we need to define models for the
attribute $\theta$ and the two deterministic processes, $\b$ and
$\v$, henceforth referred to as the mean-field processes.   
The probability mass function of the attribute $\theta$ is denoted
${\sf p}(\cdot)$. To specify the models for $\b$ and
$\v$,
we denote $\rho_t(x,u;\theta)$ as the joint distribution of the
state-action pair $(X_t,U_t)$ at time $t$, conditioned on the
attribute $\theta$.  Set
\begin{subequations}\label{eq:mf-vars}
\begin{align}
\b_t &= \sum_{\theta}  {\sf p} (\theta) \sum_{x\in\{\Aem,\Iem\}} 
\int_0^1 u 
              \rho_t(x,u;\theta) \ud u \\
\v_t &= \sum_{\theta}  {\sf p} (\theta) \big( \rho_t(\Rem;\theta) + \sum_{x\in\{\Sem,\Aem,\Iem\}}\int_0^1 u 
                   \rho_t(x,u;\theta) \ud u \big)
\end{align}
where $\rho_t(x;\theta)$ denotes the marginal, and $\beta_t$ (resp. $\alpha_t$) represent the average activity level of infected agents (resp. all agents).
\end{subequations}

\medskip
\begin{assumption}\label{ass:ass1}
Both $\b$ and $\v$ are deterministic
processes.  Furthermore, at each time $t$, $0\leq \b_t<1$ and
$0<\v_t<1$.
\end{assumption}
\medskip


\medskip

The marginal pmf $\rho_t(x;\theta)$ evolves
according to
\begin{subequations}
\begin{align}
\frac{\ud \rho_t}{\ud t} (\Sem;\theta) & = -\lambda^{\text{\tiny SA}}
                                        \b_t \int_0^1 u 
                                         \rho_t(\Sem,u;\theta) \ud u
                                         \label{eq:rho-s} \\
\frac{\ud \rho_t}{\ud t} (x;\theta) & = ({\cal
                                      A}^\dagger\rho_t)(x;\theta),\quad
                                      x\in\{\Iem,\Aem,\Rem,\Dem\} \label{eq:rho-rest}
\end{align}
\label{eq:rho}
\end{subequations}
from a given initial condition $\rho_0(x;\theta)$; 
${\cal A}^\dagger$ is the adjoint of the generator ${\cal A}$ of the
Markov process $X$.  
It is noted that control $u$ affects {\em only}
the transition from $\Sem\to\Aem$.      


\medskip
\begin{remark}
The evolution ~\eqref{eq:rho} for $\{\rho_t:t\geq 0\}$ is nonlinear because $\b_t$ 
depends upon $\rho_t$.  Apart from the terms on the righthand side
arising due to the transition $[\Sem\to\Aem]$, the other terms are
linear.  These other terms depend only upon the transition rates of the Markov
process which can themselves depend upon the attribute $\theta$.  With
$U_t\equiv1$,~\eqref{eq:rho} is an example of the classical
Kermack-McKendrick model.  The basic reproduction
number $R_0=\frac{T^{\text{r}}}{T^{\text{c}}}$ where $T^{\text{r}}$ is the
typical time until removal (i.e. recovery or death) and 
$T^{\text{c}}$ is the typical time between infectious contacts.
Evaluating $T^{\text{r}}$ and $T^{\text{c}}$ for the Markov process
with $U_t \equiv1$ and $\beta_t\equiv \bar{\b}$,
\[
R_0 = \lambda^{\text{SA}} \bar{\b}\left( \frac{\lambda^{\text{AI}}(\lambda^{\text{IR}}+\lambda^{\text{ID}})}{\lambda^{\text{AI}}+\lambda^{\text{IR}}+\lambda^{\text{ID}}}\right)^{-1}
\] 
In this paper, the choice of $U$ is
guided by an MFG formulation which is described next.
\end{remark}
\subsection{Mean-field game problem}

\noindent \textbf{Function spaces:} The filtration of the Markov
process $X$ is denoted ${\cal F}:=\{{\cal F}_t:t \geq 0\}$ where
${\cal F}_t:=\sigma(X_t)$ ($\sigma(\cdot)$ denotes the $\sigma$-algebra generated by a stochastic process).  The filtration of the observation process
$Y$ is denoted ${\cal Y} := \{ {\cal Y}_t:t\geq 0\}$ where ${\cal
  Y}_t:=\sigma(\{Y_s:0\leq s\leq t\})$.  In the two settings of the
problem, the space of admissible control inputs, denoted by ${\cal U}$, is as
follows:
\begin{align*}
\text{(fully obsvd.)} \qquad {\cal U} &= L^\infty_{\cal F}([0,\infty) ;[0,1])\\
\text{(part. obsvd.)} \qquad {\cal U} &= L^\infty_{\cal Y}([0,\infty) ;[0,1])\
\end{align*}
i.e., an admissible control input $U$ is a $[0,1]$-valued stochastic
process adapted to ${\cal F}$ in the fully observed case, and
adapted to ${\cal Y}$ in the partially observed case; the use of the
common notation ${\cal U}$ should not cause any confusion because the
two cases are treated in separate sections.  Apart from
control, the other process of interest is $(\b,\v)$ whose function
space is denoted ${\cal M}:=L^\infty ([0,\infty) ;[0,1]^2)$.

On these function spaces, define two operators:

\medskip

\begin{enumerate}
\item The operator $\Psi:{\cal M}\to {\cal U}$ as
\[
\Psi(\b,\v) = \argmin_{U\in {\cal U}} {\sf J} (U;\b,\v)
\]
\item The operator $\Xi: {\cal U} \to {\cal M}$ is 
  according to~\eqref{eq:mf-vars}-\eqref{eq:rho}.  
\end{enumerate}

Assuming the two operators are well-defined, we have:

\medskip
\begin{definition}\label{def:MFE}
A \emph{mean-field equilibrium} (MFE) is any fixed point $(\beta,\alpha)$ such that $\Xi(\Psi(\b,\v)) = (\b,\v)$.
\end{definition}
\medskip
\begin{remark}
Although $\Psi$ is well-defined under rather mild conditions, it is
difficult to justify $\Xi$ without additional assumptions on the
form of the control input $U$.  For this purpose, it is useful to note
that, provided it is well-defined, the optimal control input, denoted
$U^{\text{opt}}=\{ U_t^{\text{opt}} :t\geq 0\}$, is obtained using a deterministic feedback control law
\begin{align*}
\text{(fully obsvd.)} \qquad U_t^{\text{opt}} & = \psi_t(X_t)\\
\text{(part. obsvd.)} \qquad U_t^{\text{opt}} & = \psi_t(\pi_t)\
\end{align*}
where $\pi_t= {\sf P}(X_t\mid{\cal Y}_t)$ is the belief state; {again }the use of the
common notation $\psi_t(\cdot)$ should not cause any confusion because the
two cases are treated in separate sections. 

Now, define $\tilde{{\cal U}}\subset {\cal U}$ as the subset of all
control inputs that are obtained according to some deterministic
feedback control law ($u_t=\psi_t(x)$ or $u_t=\psi_t(\pi)$ in the two cases).  An MFE is then defined by
restricting the domain of $\Xi$ to $\tilde{{\cal U}}$.  For the fully
observed settings, such a restriction is standard.  For the partially
observed settings, we will describe an explicit construction of the
operator $\Xi$ in \Sec{sec:partially}.  
\end{remark}       

\medskip

\noindent \textbf{Notation:} The bar is used to denote
stationary (i.e., time-independent) quantities.  For example, $\psi_t
= \bar{\psi}$ means the control law is stationary, and $\b_t=\bar{\b}$
means the value of the process is a constant $\bar{\b}$.

\section{Optimality equations: Fully obsvd. case}\label{sec:fully}

\subsection{Solution for the single agent problem}

For each $x\in\Xs$ and $t\geq 0$, the value function 
\begin{align}
& v_t(x) :=  \min_{U\in L^2_{\cal F}} \nonumber \\ & {\sf E} \left( \left. 
         \int_t^{T} e^{-\gamma(s-t)}c(X_s,U_s;\v_s) \diff s
        +e^{-\gamma T} \phi(X_T) \right| X_t= x\right)
\label{eq:val_fully}
\end{align}
For $x=\Dem$ and $x=\Rem$, the value function is $v_t(\Dem)=\f(\Dem)$ and
$v_t(\Rem)=\f(\Rem)$. 
The remaining states are
$\{\Iem,\Aem,\Sem\}$.  For the state $x=\Iem$, the value function solves
the HJB equation
\begin{align*}
-\frac{\ud v_t}{\ud t} (\Iem) + (\gamma+\lambda^{\text{\tiny IR}} +
\lambda^{\text{\tiny ID}} )v_t(\Iem)  = & c^{\text{h}} (\Iem)  + \lambda^{\text{\tiny IR}}\f(\Rem)+
\lambda^{\text{\tiny ID}} \f(\Dem)  \\
& \;\; + \min_{u\in[0,1]} \{c^a(\Iem) - \v_t\} u 
\end{align*}
Now because the altruistic cost $c^a(\Iem) = 1$ and because $\v_t<1$
(Assumption~\ref{ass:ass1}), the optimal action for an infected agent
is to use $U_t^{\text{opt}}=\psi_t(\Iem)=0$, and 
\begin{align*}
-\frac{\ud v_t}{\ud t} (\Iem) + (\gamma+\lambda^{\text{\tiny IR}} +
\lambda^{\text{\tiny ID}} )v_t(\Iem)  = & c^{\text{h}} (\Iem)  + \lambda^{\text{\tiny IR}}\f(\Rem)+
\lambda^{\text{\tiny ID}} \f(\Dem)
\end{align*}
whose solution is stationary (i.e., time-independent) and given by
\[
v_t(\Iem) = \f(\Iem) := \frac{c^{\text{h}} (\Iem)}{(\gamma+\lambda^{\text{\tiny IR}} +
\lambda^{\text{\tiny ID}} )} + \frac{\lambda^{\text{\tiny IR}}\f(\Rem)+
\lambda^{\text{\tiny ID}} \f(\Dem)}{(\gamma+\lambda^{\text{\tiny IR}} +
\lambda^{\text{\tiny ID}} )}
\]
In the remainder of the paper, we make the following assumption whose
justification is provided as part of the remark after
the value function is fully described.

\medskip
\begin{assumption}\label{ass:ass2}
The value $\f(\Iem)$ is positive.  
\end{assumption}
\medskip

For the presymptomatic state $x=\Aem$, the value function solves 
\begin{align*}
-\frac{\ud v_t}{\ud t} (\Aem) + (\gamma+\lambda^{\text{\tiny AI}})v_t(\Aem) 
=c^{\text{h}}(\Aem) &+ \lambda^{\text{\tiny AI}} \f(\Iem) \\&+ \min_{u\in[0,1]} \{c^a(\Aem) - \v_t\} u
\end{align*}
and once again because of the nature on the altruistic cost,
$c^a(\Aem)=1$, the optimal action for an asymptomatic agent is to use
$U_t^{\text{opt}}=\psi_t(\Aem)=0$.  With the health cost $c^{\text{h}}(\Aem)=0$, 
\[
v_t(\Aem) = \f(\Aem):=\frac{\lambda^{\text{\tiny AI}}}{\gamma+\lambda^{\text{\tiny AI}}}\f(\Iem)
\] 
for the fully observed problem.  

It remains to obtain $v_t(\Sem)$.  The HJB equation is 
\begin{align*}
-\frac{\ud v_t}{\ud t} (\Sem) + \gamma v_t(\Sem) = \min_{u\in[0,1]}
\left( \lambda^{\text{\tiny SA}}\b_t (\f(\Aem)- v_t(\Sem))
  - \v_t\right)  u 
\label{eq:hjb-fully-s}
\end{align*}
In the following, it is assumed that a unique solution exists and
$U^{\text{opt}}_t = \psi_t(\Sem)$ is obtained as a feedback control
law.  

One may obtain additional insights by considering the stationary case
whose solution is described in the following
proposition with proof in the Appendix~\ref{appdx:pf-FO-stationary}.

\medskip

\begin{proposition}[Stationary solution] Suppose 
  $\gamma>0$, $\f(\Iem)>0$, $\b_t =
  \bar{\b}$ and $\v_t =\bar{\v}$ are both constants, 
  and $\bar{\v}<1$.  Then the optimal control for a
  susceptible agent is stationary and described by the following cases:
\begin{enumerate}
			\item If $
                          \bar{\beta}<\frac{\bar{\v}}{\lambda^{\text{\tiny
                                SA}}\f(\Aem)}$ then
			the  optimal control $U_t^{\text{opt}}=1$ and the
			optimal value $v_t(\Sem) =
                        \frac{\lambda^{\text{\tiny SA}} \bar{\beta} \f(\Aem)
                          - \bar{\v}}{\lambda^{\text{\tiny SA}}\bar{\b} + \gamma}$.  
			\item If $\bar{\b}\geq
                          \frac{\bar{\v}}{\lambda^{\text{\tiny SA}}\f(\Aem)}$ then
			the optimal control $U_t^{\text{opt}}=0$ and the optimal value $v_t(\Sem) = 0$.  
\end{enumerate}\label{prop:prop1}
\end{proposition}
\medskip
\begin{remark}
The $\gamma=0$ case is ill-posed for the following two reasons:
\begin{enumerate}
\item Suppose $\bar{\beta}<\frac{\bar{\v}}{\lambda^{\text{\tiny
        SA}}\f(\Aem)}$.  Then the optimal control $U^{\text{opt}}$ is
  not uniquely defined.  In fact, any non-zero choice of $U^{\text{opt}}$
  yields the same value $v_t(\Sem) =
                        \frac{\lambda^{\text{\tiny SA}} \bar{\beta} \f(\Aem)
                          - \bar{\v}}{\lambda^{\text{\tiny
                              SA}}\bar{\b}}$.
\item Suppose $\bar{\beta} > \frac{\bar{\v}}{\lambda^{\text{\tiny
        SA}}\f(\Aem)}$.  Then the optimal control $U_t^{\text{opt}} =
  0$ with value $v_t(\Sem) =0$.  This is directly verified from using
  the definition of the value function.  However, the HJB equation is
  not useful in this regard because, using zero control, the
  terminal time $T=\infty$; cf.,~\cite{van2007lecture}.          
\end{enumerate}

\medskip

Therefore, $\gamma$ serves as a regularization parameter.  Another
choice is to modify $\text{rate}[\Sem \to
\Aem] = \lambda^{\text{\tiny SA}} (\gamma + \b_t\,  U_t)$ which will
also serve to regularize the problem (precluding $T=\infty$ for {\em all}
choices of control).   

\end{remark}
\medskip
\begin{remark}
We next justify Assumption~\ref{ass:ass2} ($\f(\Iem)>0$).  A susceptible
agent always has an option to stay isolated (choose $U_t=0$ for all $t\geq
0$) and obtain the associated possibly sub-optimal value ${\sf J}(0) = 0$.
Assumption~\ref{ass:ass2} says that the cumulative cost of being
infected is greater than cost of staying isolated.  Without such an
assumption, an agent may wish to become active for the purposes of
getting infected and thereby lowering their value. 
Let
$\b^{\text{crit}}:=\frac{\bar{\v}}{\lambda^{\text{\tiny SA}}\f(\Aem)}$.  It
is the critical value of $\bar{\b}$ when the cost of getting infected
balances off the reward of being active, i.e., 
\[
\lambda^{\text{\tiny SA}}\b^{\text{crit}}\f(\Aem) = \bar{\v}
\]
\end{remark}






\subsection{Mean-field game}

In the fully observed version of the MFG, each agent uses
the optimal control $U^{\text{opt}}_t = \psi_t(X_t)$.  For a population with
heterogenous agents, notation $\psi_t(\cdot;\theta)$ is used to
denote the dependence on the attribute $\theta$.  Using the optimal
control, 
\[
\frac{\ud \rho_t}{\ud t} (\Sem;\theta)  = -\lambda^{\text{\tiny SA}}
                                        \b_t \psi_t(\Sem;\theta)
                                         \rho_t(\Sem;\theta) 
\]
and because the optimal control $\psi_t(x)=0$ for $x\in\{\Aem,\Iem\}$ 
\begin{align*}
\b_t &:= \sum_{\theta} {\sf p} (\theta)
              \psi_t(\Sem;\theta)  \rho_t(\Sem;\theta) \\
\v_t &:= \sum_{\theta}  {\sf p} (\theta) \left( \rho_t(\Rem;\theta) +
       \psi_t(\Sem;\theta) \rho_t(\Sem;\theta) \right)
\end{align*}
The main result is the following proposition whose straightforward
proof is omitted on account of space.

\medskip

\begin{proposition}
Suppose $\rho_0(x;\theta)$ is the initial pmf for the agents.  The 
solution for the fully observed MFG problem is:
\begin{enumerate}
\item For a single agent, the optimal control is stationary 
\[
\psi_t(x) = \bar{\psi}(x) := \begin{cases} 0 & x\in\{\Aem,\Iem\} \\
1 & x = \Sem
\end{cases}
\]
\item For the population, the distribution evolves as
\begin{align*}
\rho_t(\Sem;\theta) & = \rho_0(\Sem;\theta) \\
\frac{\ud \rho_t}{\ud t} (x;\theta) & = ({\cal
                                      A}^\dagger \rho_t)(x;\theta),\quad x\in\{\Iem,\Aem,\Rem,\Dem\}
\end{align*}

\item The consistent mean-field terms are as follows:
\begin{align*}
\b_t  = 0,\quad 
\v_t  = \sum_{\theta} {\sf p} (\theta) \left( \rho_0(\Sem;\theta) +
  \rho_t(\Rem;\theta) \right)
\label{eq:fixed-point}
\end{align*}
\end{enumerate}
\label{prop:mf-eqb}
\end{proposition}
\medskip
\begin{remark}
The conclusions of the theorem are not very practical.  It indicates 
that on the ideal planet (of vulcan) where agents are perfectly rational and have
perfect information, both the presymptomatic and symptomatic agents will
self-isolate, and therefore $R_0=0$, and a susceptible agent can continue to
party without consequence.  The main utility of the fully observed
case is to set up the problem whereby the effects of some of the
underlying assumptions -- perfect rationality and perfect information
-- can be investigated.  In the following, we consider the partially
observed problem where an agent is rational but does not have perfect
information regarding her epidemiological state.      
\end{remark}

\section{Optimality equation: Partially obsvd. case}\label{sec:partially}

\subsection{Solution for the single agent}

The partially observed problem is converted to a fully observed
one by introducing the belief state which at time $t$ is denoted by
\[
\pi_t := \begin{bmatrix} \pi_t(\Sem) & \pi_t(\Aem) &
  \pi_t(\Iem)  & \pi_t(\Rem) & \pi_t(\Dem) \end{bmatrix}
\]
where $\pi_t(x):= {\sf P} ([X_t = x] \mid {\cal Y}_t)$ for $x \in \Xs$.
Since the events $[X_t = \Iem]$ and $[X_t = \Dem]$ are both contained
in ${\cal
  Y}_t$, $\pi_t$ is not an arbitrary element of the
probability simplex in $\mathbb{R}^5$.  Let ${\cal P}^1$ denote the
set of pmf-s on $\{\Sem,\Aem\}$ and let ${\cal P}^2
= \{\delta_\Iem,\delta_\Rem,\delta_\Dem\}$.  Then the state-space for
the belief is ${\cal P}^1 \cup {\cal P}^2$. For $t\geq 0$ and 
$\mu\in {\cal P}^1 \cup {\cal P}^2$, the value function 
\begin{align*}
& v_t(\mu)  :=\min_{U\in L^2_{\cal Y}}\\
& {\sf E} \left(
           \int_t^Te^{-\gamma(s-t)} c(X_s,U_s;\v_s) {\rm
           d} s + e^{-\gamma T}\phi(X_T) \mid \pi_t = \mu  \right)
\label{eq:val_partially}
\end{align*}
There are two cases to consider: (i) when $\mu\in{\cal P}^2$; and (ii)
when $\mu\in{\cal P}^1$.  In the first case, when $\mu\in{\cal P}^2$,
the problem reduces to the fully-observed settings, and the value
function is given by
\[
v_t(\delta_\Dem) = \f(\Dem),\;\;v_t(\delta_\Rem) = \f(\Rem),\;\;v_t(\delta_\Iem) = \f(\Iem)
\] 
The optimal control for the agent in the infected state
($\pi_t=\delta_{\Iem}$) is $U_t^{\text{opt}} = \psi_t (\delta_{\Iem})
= 0$.  

For the second case, when $\mu\in{\cal P}^1$, a nonlinear filter is
used to obtain the evolution of the belief.  For this purpose,
consider first the random variable $\tau=
\tau (\omega) = \inf \{ t>0 : X_t(\omega) = \Iem \}$.  Now, $\tau$ is
a ${\cal Y}_t$-stopping time and 
\[
\pi_t= \begin{bmatrix} \pi_t(\Sem) & \pi_t(\Aem) &
  0 & 0 & 0 \end{bmatrix} \quad \text{for}\;\;t<\tau
\]
Let $A_t:=\pi_t(\Aem)$ for $t<\tau$.  Then the
stochastic process $\{A_t\in[0,1]:0\leq t<\tau\}$ evolves according to the nonlinear
filter (derived from the general form given in~\cite{confortola2013filtering}):
\[
\frac{\ud A_t}{\ud t} = (1-A_t)(\lambda^{\text{\tiny SA}} \b_t U_t -
A_t \lambda^{\text{\tiny AI}}),\;\; 0\leq t<\tau,\;\; A_0 = \pi_0(\Aem)
\]

We identify ${\cal P}^1$ with the interval $[0,1]$ with $a$ serving as its
coordinate ($a$ is the value of $A_t$).  For an arbitrary
element $\mu = [1-a,a, 0,  0, 0]$ in ${\cal P}^1$, we denote the
value function with respect to the $a$-coordinate as 
\[
\phi_t(a) := v_t(\mu), \quad 0\leq a \leq 1,\quad t\geq 0
\]
The process $\{\phi_t(a)\in\Re:0\leq a\leq 1,t\geq 0\}$ solves the HJB equation whose derivation appears in Appendix~\ref{appdx:derivations}:
\begin{align}
-\frac{\partial \phi_t}{\partial t} (a)  & + \gamma\phi_t(a)  \nonumber \\ = & - \lambda^{\text{\tiny AI}} a(1-a)
\frac{\partial \phi_t}{\partial a} (a) + \lambda^{\text{\tiny AI}} a
                                                                    (\f(\Iem)-\phi_t(a)) \nonumber \\
 & + \min_{u \in [0,1]} \left( \lambda^{\text{\tiny SA}} \b_t
   (1-a) \frac{\partial \phi_t}{\partial a}(a) + (a -\v_t)
   \right) u \label{eq:HJB_POMFG}
\end{align}
In the
following, it is assumed that a unique solution exists and yields a
well-posed optimal control law, denoted as $\psi_t(a)$ for
$a\in[0,1]$ and $t\geq 0$.  
The optimal control is
\[
U_t^{\text{opt}} = \psi_t(A_t), \quad 0\leq t<\tau
\] 
Because $A_t\equiv 1$ is an equilibrium of the filter, the HJB equation for
the point $a=1$ reduces to an ordinary differential equation whereby
\[
-\frac{\ud \phi_t}{\ud t} (1)  + \gamma\phi_t(1)  = \lambda^{\text{\tiny AI}}  (\f(\Iem)-\phi_t(1)) + \min_{u \in
  [0,1]} (1 -\bar{\v})u
\]
and because $\bar{\v}<1$ (Assumption~\ref{ass:ass2}), the optimal control law  is
$
\psi_t(1)=0
$. 
The optimal value function solves the ODE
\[
-\frac{\ud \phi_t}{\ud t} (1)  + \gamma\phi_t(1)  = \lambda^{\text{\tiny AI}}  (\f(\Iem)-\phi_t(1)) 
\]
whose solution is easily obtained as
\[
\phi_t(1) = \frac{\lambda^{\text{\tiny AI}}}{\gamma+\lambda^{\text{\tiny AI}}}\f(\Iem)
\]
It is noted that the righthand side is the value $\f(\Aem)$ for the fully
observed case.  

In summary, the optimal value is $\phi_t(1) = \f(\Aem)$
and the optimal control law is $\psi_t(1)=0$.  That is, an 
agent who has a perfect belief that she is asymptomatic will act the
same way (isolate) and will have the same value as her fully observed counterpart.  
\subsection{Partially observed mean-field game}

To setup an MFG, consider the space of probability
distributions on the belief space ${\cal P}^1\cup  {\cal P}^2$.  The
random variable $A_t$ is well-defined on the set $[t<\tau]$, and we
denote by $p_t(a)$ as its density for $0\leq a \leq 1$:
\[
{\sf P}([a< A_t <  a+\ud a]\cap [t<\tau]) = p_t(a)\ud a, \quad t\geq 0 
\]  

\medskip
\begin{assumption}\label{ass:ass3}
The density $p_t(0)=0$ for all $t\geq 0$.  
\end{assumption}
\medskip


Under Assumption~\ref{ass:ass3} that an agent uses the optimal control
$U_t=U_t^{\text{opt}} = \psi_t(A_t)$ for $0\leq t<\tau$, the density process $\{p_t(a)\in [0,\infty):0\leq
a\leq 1,t\geq 0\}$ solves the FPK equation whose derivation appears in
Appendix~\ref{appdx:derivations}:
\begin{equation}
\frac{\partial p_t}{\partial t}(a)= - \frac{\partial}{\partial
  a}\big( (1-a)(\lambda^{\text{\tiny SA}}\b_t \psi_t(a)-
  a\lambda^{\text{\tiny AI}})  p_t(a)\big)
-a\lambda^{\text{\tiny AI}}p_t(a)
\label{eq:FPK_POMFG_1}
\end{equation}
where $p_0(a)$ is the initial density (assumed given). 

By using the tower property,
\begin{align*}
\rho_t(\Aem) & = {\sf P}([X_t = \Aem]) \\
& = {\sf E}(\pi_t(\Aem)) = {\sf
  E}(A_t1_{[t<\tau]}) = \int_0^1 a p_t(a) \ud a 
\end{align*}
and therefore we have
\begin{subequations}\label{eq: FPK_POMFG_2}
\begin{align}
\frac{\ud \rho_t}{\ud t} (\Iem)&= \lambda^{\text{\tiny AI}} \int_0^1 a
                                 p_t(a)\ud a -(\lambda^{\text{\tiny
                                 ID}} + \lambda^{\text{\tiny IR}}) \rho_t(\Iem)\\
\frac{\ud \rho_t}{\ud t}(x) &= ({\cal A}^\dagger \rho_t)(x),\quad x\in\{\Rem,\Dem\}
\end{align}
\end{subequations}
where expression for ${\cal A}^\dagger$ is obtained from the
transition graph.  

With a heterogenous population, the notation $p_t(a;\theta)$ is used
to denote the
density conditioned on the attribute $\theta$ and $\psi_t(a;\theta)$
is the optimal control law.  The mean-field processes are then consistently
obtained as 
\begin{subequations}\label{eq:consistency_POMFG}
\begin{align}
\b_t & = \sum_{\theta} {\sf p}(\theta) \int_0^1 a \psi_t(a;\theta) p_t(a;\theta) \ud a
  \\
\v_t & = \sum_{\theta} {\sf p}(\theta)  \left( \rho_t(\Rem;\theta) + \int_0^1
       \psi_t(a;\theta) p_t(a;\theta) \ud a \right)
\end{align}
\end{subequations}

This completes the derivation of the system of equations for the
partially observed MFG: Eq.~\eqref{eq:FPK_POMFG_1}-\eqref{eq: FPK_POMFG_2} is the forward FPK equation.
Eq.~\eqref{eq:HJB_POMFG} is the backward HJB equation.
Eq.~\eqref{eq:consistency_POMFG} defines the consistency relationship
that links the two equations.  Its solution is an MFE (satisfies Defn.~\ref{def:MFE}).  

The analytical and numerical study of the mean-field equations is a
subject of continuing work. In the following, we discuss some
preliminary analytical results.


\subsection{Some special cases}

\def\s{\bar{\psi}}
To gain further insights into the model, we consider stationary
solutions of the HJB equation~\eqref{eq:HJB_POMFG}.  For this purpose, in this subsection,
we assume that $\b_t = \bar{\b}$ and $\v_t =\bar{\v}$ are both constants.
The stationary HJB equation is then
\begin{align}
\gamma\f(a)  = & - \lambda^{\text{\tiny AI}} a(1-a)
\frac{\ud \f}{\ud a} (a) + \lambda^{\text{\tiny AI}} a (\f(\Iem)-\f(a)) \nonumber\\
 & + \min_{u \in [0,1]} \left( \lambda^{\text{\tiny SA}} \bar{\b}
   (1-a) \frac{\ud \f}{\ud a}(a) + (a -\bar{\v})
   \right) u 
\label{eq:stationary-HJB}
\end{align}
The stationary optimal control law, obtained upon evaluating the minimizer of the HJB equation,
is denoted by $u = \s(a)$ for $a\in[0,1]$.  As already described for the
general non-stationary case, at $a=1$, the optimal value is 
$\f(1)=\frac{\lambda^{\text{\tiny AI}}}{\gamma+\lambda^{\text{\tiny AI}}}\f(\Iem)$
and the associated optimal control is $\s(1) = 0$.  

\medskip

\newP{Case 1. Limit $\lambda^{\text{\tiny AI}} \uparrow \infty$}  Evaluate
the stationary HJB equation at $a=0$:
\[
\gamma \f(0) = \min_{u\in[0,1]}\left(\lambda^{\text{\tiny SA}}\bar{\b}\frac{\ud \f}{\ud a}(0)-\bar{\v}\right)u
\] 
Assume now that $\phi$ is continuously differentiable at $a=0$.  Then
the value of $\frac{\ud \f}{\ud a}(0)$ can be obtained using the
dominant balance in the limit as $\lambda^{\text{\tiny AI}} \to \infty$ and
$a>0$ but small:
\[
- \lambda^{\text{\tiny AI}} a(1-a)
\frac{\ud \f}{\ud a} (a) + \lambda^{\text{\tiny AI}} a (\f(\Iem)-\f(a)) = 0
\]
which yields $\frac{\ud \f}{\ud
  a}(0)=\f(\Iem)-\f(0)$.  Therefore, in the limit as
$\lambda^{\text{\tiny AI}} \uparrow \infty$, the stationary HJB equation for
$a=0$ is given by
\[
\gamma \f(0) =
\min_{u\in[0,1]}\left(\lambda^{\text{\tiny SA}}\bar{\b}(\f(\Iem)-\f(0))-\bar{\v}\right)
u
\]
which is identical to the stationary HJB
equation~\eqref{eq:stat_HJB_fully_obs_phi} for the fully observed
problem.  The optimal control  law is
\[
\s(0) = \begin{cases} 
1 & \text{if}\;\; \lambda^{\text{\tiny SA}}\bar{\b}\f(\Iem) < \bar{\v}\\
0 & \text{o.w.}
\end{cases}
\]
which is also the same as the optimal control for a susceptible agent in
the fully observed case.  This suggests that if $a= 0$ (agent
has certain belief that she is susceptible), she will be active if
$\bar{\b} < \bar{\b}^{\text{crit}}$.  

To obtain insights for arbitrary values of $a\neq 1$, consider
the nonlinear filter using the stationary control law $\s$:
\begin{align}
\frac{\ud A_t}{\ud t} = (1-A_t)(\lambda^{\text{\tiny SA}} \bar{\b} \s(A_t) -
A_t \lambda^{\text{\tiny AI}}),\;\; 0< t<\tau,\;\; A_0 \neq 1
\label{eq:filter}
\end{align}
where we note that $0\leq \s(A_t) \leq 1$ and $\s(1)=0$.  In the asymptotic
limit as $\lambda^{\text{\tiny AI}} \to \infty$, its solution is given by
\[
A_t = \frac{\lambda^{\text{\tiny SA}}\bar{\b}}{\lambda^{\text{\tiny AI}}} \s(0) +
o ( \frac{1}{\lambda^{\text{\tiny AI}}})\approx 0,\quad 0< t<\tau
\]
and therefore, the optimal control $U_t^{\text{opt}} = \s(A_t) \approx
\s(0)$ whenever $A_0\neq 1$.  

In summary, an agent who has a perfect belief that she is asymptomatic
$(a=1)$ will isolate ($\s(1)=0$).  For all other values ($a\neq 1$),
the agent will behave as a susceptible agent in the fully
observed settings of the problem:  isolate if $\bar{\b} >
\bar{\b}^{\text{crit}}$ and fully active if $\bar{\b} <
\bar{\b}^{\text{crit}}$.




\medskip

\newP{Case 2. $\lambda^{\text{\tiny AI}}$ large} The limit
($\lambda^{\text{\tiny AI}} = \infty$) represents the case when there is no
uncertainty in belief.  
As one deviates away from the
limit, the uncertainty increases and an agent is no longer perfectly 
sure of her epidemiological state.  The following proposition
shows that for sufficiently large values of the parameter
$\lambda^{\text{\tiny AI}}$, the stationary optimal control law is of threshold type.
\medskip
\begin{proposition}\label{prop:partial_threshold}
Suppose $\b_t =
  \bar{\b}$, $\v_t =\bar{\v}$.  Then
\begin{enumerate}
\item If $\bar{\b}\geq \bar{\b}^{\text{crit}}$
  then the optimal control law is ${\s}(a) = 0$ for all $0\leq a \leq 1$.
\item For each fixed  $\bar{\beta}<\frac{\bar{\v}}{\lambda^{\text{\tiny
        SA}}\f(\Iem)}<
  \bar{\b}^{\text{crit}}$ there exists a
  $\underline{\lambda}^{\text{\tiny AI}} = \underline{\lambda}^{\text{\tiny AI}}
  (\bar{\b})$ such that for all $\lambda^{\text{\tiny AI}} >
  \underline{\lambda}^{\text{\tiny AI}}$, the optimal control law is of
  threshold type:
\begin{align}
\s(a)= \begin{cases}
1 & \text{if} \;\; 0\leq a< {a}^{\text{thresh}} \\
0 & \text{if} \;\; {a}^{\text{thresh}} < a \leq 1
\end{cases}
\label{eq:sta_policy}
\end{align}
where the threshold $ {a}^{\text{thresh}}\in(0,1)$.  
The function $\underline{\lambda}^{\text{\tiny AI}}
  (\bar{\b})$ is monotonic in its argument and $\lim_{\bar{\b}\downarrow 0}
    \underline{\lambda}^{\text{\tiny AI}}
  (\bar{\b})=0$.  
  
\end{enumerate}
\end{proposition}

\medskip

The exact formulae for the function $\underline{\lambda}^{\text{\tiny AI}}
  (\bar{\b})$ and the threshold $ {a}^{\text{thresh}}$ are fairly
  complicated.  These formulae appear along with the proof of
  Prop.~\ref{prop:partial_threshold} in the
  Appendix~\ref{appdx:pf-FO-MFG}.
Apart from these special cases, we do not yet have a complete
understanding of the stationary solutions
of~\eqref{eq:stationary-HJB}.  This remains a subject of continuing
work.  
\section{Conclusions and directions for future work}

In this paper, we proposed a partially observed MFG model
for epidemics.  The main contribution is derivation of the
forward-backward equations.  The analytical and numerical study of
these equations is 
a topic of continuing research.  
A major simplifying assumption in the model is that we ignored the
transition from $\Aem$ to $\Rem$.  With the transition present, an agent will
maintain a belief over a three-dimensional state ($\Sem,\Aem,\Rem$)
leading to a loss of total order on beliefs.  Although the extension
of the forward-backward equations should be easily possible, the
analysis will be much more complicated.  
Including the effects of
testing and vaccination are other directions to extend the basic model.



  
\bibliographystyle{IEEEtran}
\bibliography{references,referencesnew}

\appendix

	\subsection{Proof of
          \Prop{prop:prop1}}\label{appdx:pf-FO-stationary}

With $\b_t = \bar{\b}$ and $\v_t = \bar{\v}$, the HJB equation for $v_t(\Sem)$ is 
\[
-\frac{\ud v_t}{\ud t} (\Sem) + \gamma v_t(\Sem) = \min_{u\in[0,1]}
\left( \lambda^{\text{\tiny SA}}\bar{\b}(\f(\Aem)- v_t(\Sem))
  - \bar{\v} \right)  u 
\]
We investigate its stationary solutions $v_t(\Sem) =
\bar{\phi}(\Sem)$ in which case the stationary HJB equation is
\begin{equation}\label{eq:stat_HJB_fully_obs_phi}
\gamma \bar{\phi}(\Sem)= \min_{u\in[0,1]}
\underbrace{\left( \lambda^{\text{\tiny SA}}\bar{\b}(\f(\Aem)- \bar{\phi}(\Sem))
  - \bar{\v} \right)}_{=:M}  u 
\end{equation}
We have the  following two cases:
\begin{itemize}
\item If $\bar{\beta}<\frac{\bar{\v}}{\lambda^{\text{\tiny
        SA}}\f(\Aem)}$ then $\bar{\phi}(\Sem) =
  \frac{\lambda^{\text{\tiny SA}} \bar{\beta} \f(\Aem) -
    \bar{\v}}{\lambda^{\text{\tiny SA}}\bar{\b} + \gamma}$ solves the HJB equation with the minimizing choice of
$u=1$ because
	\[
	M
        = \frac{\gamma}{\lambda^{\text{\tiny SA}}\bar{\b} +
          \gamma}\big(\lambda^{\text{\tiny SA}} \bar{\beta} \f(\Aem) - \bar{\v}\big)<0
	\]
\item If $\bar{\beta} > \frac{\bar{\v}}{\lambda^{\text{\tiny
        SA}}\f(\Aem)}$ then $\bar{\phi}(\Sem) =0$  solves the HJB equation with the minimizing choice of
$u=0$ because
\[
M
        = \lambda^{\text{\tiny SA}}\bar{\b}\f(\Aem)-\bar{\v} > 0
\]
\end{itemize}
\subsection{Derivations}\label{appdx:derivations}
Let ${\cal P}:={\cal P}^1\cup {\cal P}^2$.  
The belief process $\pi:=\{\pi_t\in{\cal P}:t\geq
0\}$ is a Markov process~\cite[Theorem 1.7]{xiong2008book}.  The HJB equations and FPK equations
are easily derived once we obtain the infinitesimal generator of the
process.  

\medskip

\newP{Infinitesimal generator} Consider a smooth test function $v:{\cal P} \to \Re$. 
Let $\mu\in {\cal P} $.  The infinitesimal generator
(for the general time-inhomogeneous case) is
	\[
	({\cal A}_t^u v)(\mu) = \lim_{\delta t\downarrow 0}
        \frac{\E(v(\pi_{t+\delta t})|\pi_t=\mu)-v(\mu)}{h}
	\]
There are two cases to consider:

\medskip

\noindent $\bullet$ If $\mu\in {\cal P}^2$, upon identifying the
measures $\{\delta_\Iem,\delta_\Rem,\delta_\Dem\}$ with the states
$\{\Iem,\Rem,\Dem\}$, the generator is the same as the generator for
the Markov process.  

\medskip

\noindent $\bullet$
If $\mu\in {\cal P}^1$ then using the coordinate $a$ for ${\cal
    P}^1$, $\mu = [1-a,a,0,0,0]$ for $a\in[0,1]$.  With $\pi_t=\mu$,
  in the asymptotic limit as $\delta t\to 0$,
\[
\pi_{t+\delta t} = \begin{cases}
	\mu + e f(a,t) \delta t + o(\delta t) & \text{w.p.} \;\;
        (1-a\lambda^{\text{\tiny AI}}\delta t) + o(\delta t)\\
	\delta_\Iem &\text{w.p.} \;\;
        a\lambda^{\text{\tiny AI}}\delta t + o(\delta t)
	\end{cases}
	\]   
where $f(a,t)=(1-a)(\lambda^{\text{\tiny SA}} \beta_t u-a
\lambda^{\text{\tiny AI}})$ and $e = [-1,1,0,0,0]$. 
Denoting $v(\mu) = \phi(a)$, the generator is then easily calculated to be
\[
({\cal A}_t^u v)(\mu) = f(a,t)  \frac{\partial \phi}{\partial a}(a) +
a \lambda^{\text{\tiny AI}} (v(\delta_\Iem)-\phi(a) )
\]
where the superscript $u$ denotes the fact that $f(a,t)$, and therefore 
also the generator, depends also upon
$u$.  The subscript $t$ denotes the fact that the generator is for a
time-inhomogeneous Markov process (because $\beta$ may depend
upon time).  

\medskip

\newP{Derivation of the HJB equation} 
The HJB equation is
\[
- \frac{\partial v_t}{\partial t} (\mu) + \gamma v_t(\mu) = \min_{u\in[0,1]} \left( {\cal A}_t^u v_t(\mu) + \mu
  (c(\cdot,u;\v_t)) \right)
\]
For $\mu = [1-a,a,0,0,0] $ and $v_t(\mu) = \phi_t(a)$, the HJB
equation~\eqref{eq:HJB_POMFG} is obtained because 
\[
\mu\big(c(\cdot,u;\v_t)\big) = (1-a)(-\v_t u) + a (1-\v_t)u = (a - \v_t) u
\]

\medskip

\newP{Derivation of the FPK equation}
We derive the adjoint of the generator ${\cal A}^u$ where dependence
on $t$ is suppressed for notational ease. Let $\rho$ be a measure on ${\cal P}$. On ${\cal P}^1$, $\rho$ has density $p(a)$. Consider $\rho({\cal A}^uv) = \int {\cal A}^uv(\mu) \rho(\ud \mu)=$ 
\begin{align*}
&\int_0^1 p(a)\bigg[f(a)\frac{\partial \phi}{\partial a} +a\lambda^{\text{\tiny AI}}\big(v(\delta_\Iem)-\phi(a)\big) \bigg]\ud a\\
&\quad  -\rho(\delta_\Iem)(\lambda^{\text{\tiny IR}}+\lambda^{\text{\tiny ID}})v(\delta_\Iem) + \rho(\delta_\Rem)\lambda^{\text{\tiny IR}}v(\delta_\Rem) + \rho(\delta_\Dem)\lambda^{\text{\tiny ID}}v(\delta_\Dem)\\
= & - \int_0^1\phi(a)\left(
  \frac{\partial}{\partial a}\big(p(a)f(a) \big)
  +a\lambda^{\text{\tiny AI}}p(a)\right) \ud a \\
& +  p(a)\phi(a)f(a)\Big\vert_0^1  +
  v(\delta_\Rem)\lambda^{\text{\tiny IR}}\rho(\delta_\Rem) +
  v(\delta_\Dem)\lambda^{\text{\tiny ID}}\rho(\delta_\Dem) \\
& + v(\delta_\Iem)\left(\int_0^1a\lambda^{\text{\tiny AI}}p(a)\ud a-\rho(\delta_\Iem)(\lambda^{\text{\tiny IR}}+\lambda^{\text{\tiny ID}})\right)
\end{align*}
where boundary terms vanish because $f(1)=0$ and $p(0)=0$
(Assumption~\ref{ass:ass3}).


\subsection{Proof of \Prop{prop:partial_threshold}}\label{appdx:pf-FO-MFG}

We are interested in solutions of~\eqref{eq:stationary-HJB} repeated below
\begin{align}
\gamma\f(a)  = & - \lambda^{\text{\tiny AI}} a(1-a)
\frac{\ud \f}{\ud a} (a) + \lambda^{\text{\tiny AI}} a (\f(\Iem)-\f(a)) \nonumber\\
 & + \min_{u \in [0,1]} \underbrace{\left( \lambda^{\text{\tiny SA}} \bar{\b}
   (1-a) \frac{\ud \f}{\ud a}(a) + (a -\bar{\v})
   \right)}_{:= \mathcal{M} (a)} u
\label{eq:stationary-HJB-rep}
\end{align}

\medskip

\noindent\textbf{Proof of Part 1):} If $\bar{\b}\geq
\bar{\b}^{\text{crit}}$ then $\f(a) = \f(\Aem)a$
solves~\eqref{eq:stationary-HJB-rep} with $u = 0$.  To show that
$u=0$ is a minimizer, we need to verify that $\mathcal{M}(a)>0$ for
all $0\leq a \leq 1$.  This
is true because $$\mathcal{M}(a)= \lambda^{\text{\tiny SA}} \bar{\b} (1-a)
\f(\Aem) +  (a -\bar{\v})$$ is an affine function with $\mathcal{M}(1)
= 1 - \bar{\v} > 0$ (Assumption~\ref{ass:ass1}) and $\mathcal{M}(0)
= \lambda^{\text{\tiny SA}} \bar{\b} \f(\Aem) - \bar{\v} > 0$ (because
$\bar{\b}\geq
\bar{\b}^{\text{crit}}$).

	
\medskip

\noindent \textbf{Proof of Part 2):} With $U_t \equiv 1$ the point $\bar{a}
:= \frac{\lambda^{\text{\tiny SA}}\bar{\b}}{\lambda^{\text{\tiny
      AI}}}$ is a stable equilibrium of the
filter~\eqref{eq:filter}. Now let 
\[
{a}^{\text{thresh}}:=\frac{\bar{\v}-\lambda^{\text{\tiny
      SA}}\bar{\b}k}{1-\lambda^{\text{\tiny SA}}\bar{\b}k}
\] 
where $k := \frac{\lambda^{\text{\tiny AI}}\f(\Iem)+1
  -\lambda^{\text{\tiny AI}}y}{\lambda^{\text{\tiny
      SA}}\bar{\b}+\gamma}$ and $y := \frac{\lambda^{\text{\tiny
      AI}}\f(\Iem)+1-\bar{\v}}{\lambda^{\text{\tiny AI}}+\gamma}$.
      
      We present the proof in three steps. In step 1, we show that if $0<\bar{a}<{a}^{\text{thresh}}<1$, then the optimal control obtained as a solution to the HJB equation is of
		threshold type~\eqref{eq:sta_policy}. In step 2, we present a tight bound to attain $0<\bar{a}<{a}^{\text{thresh}}<1$, which is implied by the inequality $\lambda^{\text{\tiny AI}} > \underline{\lambda}^{\text{\tiny AI}}(\bar{\b})$ derived in the final step 3.

\medskip

\noindent \textbf{Step 1:} Assume $0<\bar{a}<{a}^{\text{thresh}}<1$. We want to show that the solution to~\eqref{eq:stationary-HJB-rep} is the following:
\[
\f(a) = \begin{cases}
k(a-1)+y\quad &\text {if }a < {a}^{\text{thresh}}\\
\f(\Aem)a + c(1-a)^{1+b}a^{-b}\quad &\text {if }a \ge {a}^{\text{thresh}}
\end{cases}
\]
where $b:=\frac{\gamma}{\lambda^\text{\tiny AI}}$ and
$$c := \frac{\f(\Aem)-k}{(1-a^{\text{thresh}})^b (a^{\text{thresh}})^{-b-1} (b+a^{thresh})}$$
and the optimal policy is~\eqref{eq:sta_policy}.

For all $a < {a}^{\text{thresh}}$, we have
\[
(1-a)\lambda^{\text{\tiny SA}}\bar{\b}k + a-\bar{\v} < 0
\]
and then $\f(a) = k(a-1)+y$ solves the HJB equation with $u=1$. Note that this affine formula uniquely solves~\eqref{eq:stationary-HJB-rep} with $u=1$, since the homogeneous solution blows up at $\bar{a}$.

With $u=0$, we already know that $\f(\Aem)a$ solves~\eqref{eq:stationary-HJB-rep} and $c(1-a)^{1+b}a^{-b}$ is the homogeneous solution for~\eqref{eq:stationary-HJB-rep} that makes $\mathcal{M}(a)$ (and $\f(a)$)  continuous at $a^{\text{thresh}}$. Next we want to show that $\mathcal{M}(a)>0$ for all $a > a^{\text{thresh}}$. For $a > a^{\text{thresh}}$, $\mathcal{M}(a)$ can be expressed as
\begin{align*}
\mathcal{M}(a) &= \lambda^{\text{\tiny SA}} \bar{\b}
   (1-a)\left( \f(\Aem) - c (1-a)^b a^{-b} \left( 1+\frac{b}{a} \right) \right) \\
   &+ (a -\bar{\v})
\end{align*}
We readily observe that
\begin{align*}
\frac{d^2 \mathcal{M}}{da^2} (a) = \lambda^{\text{\tiny SA}} \bar{\b} \left| c \right| \left[ 1+b+(b+2) \frac{1-a}{a}  \right] \frac{d^2 \f^h}{da^2} (a) > 0
\end{align*}
where $\f^h(a):= (1-a)^{1+b}a^{-b}$. Noting that $\mathcal{M}(a^{\text{thresh}}) = 0$ and $\mathcal{M}(a)$ is convex, it remains to show that $\lim_{a \downarrow a^{\text{thresh}}} \frac{\ud \mathcal{M}}{\ud a} (a) \geq 0$.
\begin{align*}
\lim_{a \downarrow a^{\text{thresh}}} \frac{\ud \mathcal{M}}{\ud a} (a) &= 1-\lambda^{\text{\tiny SA}} \bar{\b}k - \lambda^{\text{\tiny SA}} \bar{\b} \frac{k-\f(\Aem)}{a^{\text{thresh}} (b+a^{\text{thresh}})}b(b+1)\\
&>1-\lambda^{\text{\tiny SA}} \bar{\b}k -\frac{k-\f(\Aem)}{b+a^{\text{thresh}}} \gamma (b+1)
\end{align*}
where the last inequality follows from $a^{\text{thresh}} > \bar{a}$ which implies $-b>-\frac{\gamma a^{\text{thresh}}}{\lambda^{\text{\tiny SA}} \bar{\b}}$. Therefore it remains to show
\begin{align*}
(b+a^{\text{thresh}})(1-\lambda^{\text{\tiny SA}} \bar{\b}k) -(k-\f(\Aem))\gamma (b+1) \geq 0
\end{align*}
The lefthand side of the equation above turns out to be exactly $0$. This can be shown as
\begin{align*}
&(b+a^{\text{thresh}})(1-\lambda^{\text{\tiny SA}} \bar{\b}k) -(k-\f(\Aem))\gamma (b+1)\\
&=b (1-\lambda^{\text{\tiny SA}} \bar{\b}k) + (\bar{\v}-\lambda^{\text{\tiny SA}} \bar{\b}k) - (k-\f(\Aem)) \gamma (b+1)\\
&=\bar{\v}+b-(b+1)\left((\lambda^{\text{\tiny SA}} \bar{\b}+\gamma)k-\gamma \f(\Aem)\right)\\
&=\bar{\v}+b-(b+1)\left( 1 - \frac{1-\bar{\v}}{b+1} \right)\\
&= \bar{\v}+b-b-1+1-\bar{\v} = 0
\end{align*}
and thus $\mathcal{M}(a)>0$ for all $a>a^{\text{thresh}}$.
\medskip

\noindent \textbf{Step 2:} We now show if $\lambda^{\text{\tiny AI}}>\lambda^{\text{\tiny SA}}\bar{\b}$, and
\begin{align}
\frac{1-\bar{\v}}{\lambda^{\text{\tiny AI}}+\gamma}+\frac{\lambda^{\text{\tiny SA}}\bar{\b}+\gamma}{(\lambda^{\text{\tiny AI}}-\lambda^{\text{\tiny SA}}\bar{\b})\gamma}< \frac{\bar{\v}}{\lambda^{\text{\tiny SA}}\bar{\b}} -\frac{\lambda^{\text{\tiny AI}}\f(\Iem)}{\lambda^{\text{\tiny AI}}+\gamma}
\label{eq:step2_key}
\end{align}
then $0<\bar{a}<{a}^{\text{thresh}}<1$.

Equation~\eqref{eq:step2_key} could be rearranged as
\begin{align} \label{eq:rearranged}
\frac{\lambda^{\text{\tiny AI}}\phi(\Iem)+1 -\bar{\v}}{\lambda^{\text{\tiny AI}}+\gamma} < \frac{\bar{\v}}{\lambda^{\text{\tiny SA}}\bar{\b}} - \frac{\lambda^{\text{\tiny SA}}\bar{\b} + \gamma}{(\lambda^{\text{\tiny AI}}-\lambda^{\text{\tiny SA}}\bar{\b})\gamma}(1-\bar{\v})
\end{align}
Writing $k$ as
\begin{align*}
k = \frac{1}{\lambda^{\text{\tiny SA}} \bar{\b} + \gamma} \left( \gamma \frac{\lambda^{\text{\tiny AI}}\phi(\Iem)+1 -\bar{\v}}{\lambda^{\text{\tiny AI}}+\gamma} + \bar{\v} \right)
\end{align*}
and applying~\eqref{eq:rearranged} gives
\begin{align*}
\lambda^{\text{\tiny SA}} \bar{\b} k < \bar{\v}-\frac{\lambda^{\text{\tiny SA}} \bar{\b}}{\lambda^{\text{\tiny AI}}-\lambda^{\text{\tiny SA}} \bar{\b}}
\end{align*}
Since the second term is positive we obtain $\bar{\v}>\lambda^{\text{\tiny SA}} \bar{\b} k$ which implies $0<a^{\text{thresh}}<1$. Moreover, the former inequality implies
$$(\lambda^{\text{\tiny AI}}-\lambda^{\text{\tiny SA}}\bar{\b})\lambda^{\text{\tiny SA}}\bar{\b}k < \lambda^{\text{\tiny AI}}\bar{\v}-\lambda^{\text{\tiny SA}}\bar{\b}\alpha$$
which gives
$
{\bar{a}} < {a}^{\text{thresh}}.
$

\medskip

\noindent \textbf{Step 3:} It remains to derive a sufficient condition that $\lambda^{\text{\tiny AI}} > \underline{\lambda}^{\text{\tiny AI}}(\bar{\b})$ implies~\eqref{eq:step2_key}
where
\begin{align*}
\underline{\lambda}^{\text{\tiny AI}} (\bar{\b})= \max \left(\lambda^{\text{\tiny SA}}\bar{\b}, \frac{\lambda^{\text{\tiny SA}}\bar{\b}\gamma}{(2-\bar{\v})\gamma+ \lambda^{\text{\tiny SA}}\bar{\b}} + \frac{\lambda^{\text{\tiny SA}}\bar{\b}\gamma}{\bar{\v}- \lambda^{\text{\tiny SA}}\bar{\b}\f(\Iem)}\right)
\end{align*}
Note that
\[
\frac{1-\bar{\v}}{\lambda^{\text{\tiny AI}}+ \gamma}\frac{{\lambda^{\text{\tiny AI}}+ \gamma}}{\lambda^{\text{\tiny AI}}-\lambda^{\text{\tiny SA}}\bar{\b}} + \frac{\lambda^{\text{\tiny SA}}\bar{\b}+ \gamma}{\gamma(\lambda^{\text{\tiny AI}}-\lambda^{\text{\tiny SA}}\bar{\b})} < \frac{\bar{\v}}{\lambda^{\text{\tiny SA}}\bar{\b}}-\f(\Iem)
\]
implies~\eqref{eq:step2_key} and doing simple manipulations gives
\[
\lambda^{\text{\tiny AI}} > \frac{\lambda^{\text{\tiny SA}}\bar{\b}\gamma}{(2-\bar{\v})\gamma+ \lambda^{\text{\tiny SA}}\bar{\b}} + \frac{\lambda^{\text{\tiny SA}}\bar{\b}\gamma}{\bar{\v}- \lambda^{\text{\tiny SA}}\bar{\b}\f(\Iem)}
\]
Therefore $\lambda^{\text{\tiny AI}} > \underline{\lambda}^{\text{\tiny AI}}(\bar{\b})$ implies~\eqref{eq:step2_key}. Note that the function $\underline{\lambda}^{\text{\tiny AI}}(\bar{\b})$ is monotonic in its argument and $\lim_{\bar{\b}\downarrow 0}
			\underline{\lambda}^{\text{\tiny AI}}
			(\bar{\b})=0$.

\end{document}